\documentclass[oneside, 12pt]{amsart}
\usepackage{amscd, amssymb, amsmath, mathrsfs}
\usepackage[english]{babel}
\usepackage{booktabs}
\usepackage{tikz-cd}
\usepackage{url}
\usepackage[pdftex, colorlinks=true,  citecolor=blue, linkcolor=blue, linktocpage=true]{hyperref}

\newcommand\mylabel[1]{\label{#1}\marginpar{\vspace{-1ex}\medskip\medskip\footnotesize \tt #1}}
\renewcommand\mylabel[1]{\label{#1}}
\newcommand{\mydate}{
\number\day\space
\ifcase\month \or January\or February\or March\or April\or May\or June\or July\or August\or September\or October\or November\or December\fi 
\space\number\year}

\DeclareUrlCommand\arXiv{\urlstyle{same}}

\newtheorem{theorem}{Theorem}[section]
\newtheorem*{maintheorem}{Theorem}
\newtheorem{lemma}[theorem]{Lemma}
\newtheorem{proposition}[theorem]{Proposition}
\newtheorem{corollary}[theorem]{Corollary}

\theoremstyle{definition}

\newtheorem*{acknowledgement}{Acknowledgement}

\theoremstyle{remark}

\newcommand{\ZZ}{\mathbb{Z}}

\newcommand{\CC}{\mathbb{C}}

\newcommand{\PP}{\mathbb{P}}

\newcommand{\GG}{\mathbb{G}}

\newcommand{\shE}{\mathscr{E}}
\newcommand{\shF}{\mathscr{F}}

\newcommand{\shN}{\mathscr{N}}
\newcommand{\shL}{\mathscr{L}}

\newcommand{\Aut}{\operatorname{Aut}}

\newcommand{\Bl}{\operatorname{Bl}}

\newcommand{\CH}{\operatorname{CH}}

\newcommand{\depth}{\operatorname{depth}}

\newcommand{\Ext}{\operatorname{Ext}}

\newcommand{\Gal}{\operatorname{Gal}}

\newcommand{\id}{{\operatorname{id}}}
\newcommand{\Image}{\operatorname{Im}}

\newcommand{\Kernel}{\operatorname{Ker}}

\newcommand{\lra}{\longrightarrow}

\newcommand{\maxid}{\mathfrak{m}}

\renewcommand{\O}{\mathscr{O}}

\newcommand{\Pic}{\operatorname{Pic}}

\newcommand{\pmap}{{[p]}}

\newcommand{\pr}{\operatorname{pr}}
\newcommand{\Proj}{\operatorname{Proj}}

\newcommand{\quadand}{\quad\text{and}\quad}

\newcommand{\ra}{\rightarrow}

\newcommand{\Reg}{\operatorname{Reg}}

\newcommand{\sep}{{\operatorname{sep}}}

\newcommand{\Sing}{\operatorname{Sing}}

\newcommand{\Spec}{\operatorname{Spec}}

\newcommand{\Sym}{\operatorname{Sym}}

\newcommand{\uH}{\underline{H}}

\newcommand{\uHom}{\underline{\operatorname{Hom}}}

\begin{document}

\title[Sign involutions]
      {Sign involutions on para-abelian varieties}.

\author[]{Jakob Bergqvist}
\address{Heinrich Heine University Düsseldorf, Faculty of Mathematics and Natural Sciences, Mathematical Institute, 40204 D\"usseldorf, Germany}
\curraddr{}
\email{Jakob.Bergqvist@hhu.de}

\author[]{Thuong  Dang}
\address{Heinrich Heine University Düsseldorf, Faculty of Mathematics and Natural Sciences, Mathematical Institute, 40204 D\"usseldorf, Germany}
\curraddr{}
\email{dangt@uni-duesseldorf.de}

\author[Stefan Schr\"oer]{Stefan Schr\"oer}
\address{Heinrich Heine University Düsseldorf, Faculty of Mathematics and Natural Sciences, Mathematical Institute, 40204 D\"usseldorf, Germany}
\curraddr{}
\email{schroeer@math.uni-duesseldorf.de}

% 14L30 Group actions on varieties or schemes (quotients)
% 14K15 Arithmetic ground fields for abelian varieties
% 14K30 Picard schemes, higher Jacobians
% 14J26 Rational and ruled surfaces

%\dedicatory{Preliminary version, \mydate}
\dedicatory{Second revised version, 3 April 2024}

\begin{abstract}
We study the  so-called sign  involutions on twisted forms of abelian varieties,
and show that such a sign involution exists if and only if the   class in the Weil--Ch\^{a}telet group is annihilated by two.
If these equivalent conditions hold, we prove that the Picard scheme of the quotient is \'etale and contains
no points of finite order. In dimension one, such    quotients are Brauer--Severi curves,
and we analyze the ensuing embeddings of the genus-one curve into twisted forms of Hirzebruch surfaces and weighted projective spaces.
\end{abstract}

\maketitle
\tableofcontents

\newcommand{\Inv}{\operatorname{Inv}}
\newcommand{\ind}{\operatorname{ind}}
\newcommand{\Bs}{\operatorname{Bs}}
\newcommand{\per}{\operatorname{per}}

%===========================================================
\section*{Introduction}
\mylabel{introduction}

Recall that an \emph{abelian variety} $A$ over a ground field $k$ is a group scheme that is proper, smooth,
and connected. As  a non-trivial consequence, the group law is  commutative,
such that $A$ comes with a canonical automorphism $x\mapsto -x$, the \emph{sign involution}.
Note that over the field $k=\CC$ of complex numbers, the abelian varieties correspond to complex tori $\CC^g/\Lambda$, where $\Lambda$ is a full
lattice admitting a polarization. An excellent exposition of the theory was given by Mumford \cite{Mumford 1970}.

Abelian varieties play a fundamental role in algebraic geometry, since they are basic building blocks for
algebraic groups.
In particular, for every proper scheme $X$ the  Picard group, viewed as a group scheme,
 contains a maximal abelian subvariety $A=\Pic^\alpha_{X/k}\subset\Pic_{X/k}$,
which encodes crucial geometric information (\cite{Laurent; Schroeer 2024}, Section 7 and \cite{Schroeer 2024}, Section 3). For smooth curves $X$,
these are the \emph{jacobian varieties}.
Abelian varieties are also important objects in arithmetic geometry, where the ground field  could be  a number field or a function field.
Geometric and arithmetic aspects are strongly interrelated: In  fibrations $f:Y\ra B$ of proper schemes, one     has to understand the generic fiber $X=f^{-1}(\eta)$ as a scheme
over the function field $k(B)$ of the base.

The  sign involution  $\sigma(x)=-x$ on   abelian varieties $A$ plays an important role, 
because it gives rise to the   notion of  symmetric   sheaves.
Furthermore, one can form the quotient $A/G$ for the corresponding
group $G=\{\pm 1\}$ of order two. In dimension $g=1$ this gives  the projective line,
whereas for  $g=2$ we get   Kummer surfaces, a fascinating topic going back to the 19th century.
In characteristic  $p\neq 2$ Kummer surfaces are  K3 surfaces with rational double points.
The case $p=2$ requires  extra attention, because than $A/G$ may also be
a rational surface with an elliptic singularity (\cite{Shioda 1974} and \cite{Katsura 1978}).
This is a prime example of a \emph{wild quotient singularity} (see for example
\cite{Lorenzini; Schroeer 2020} and \cite{Lorenzini; Schroeer 2023} for more on this topic). To our best knowledge,
no resolution of singularities is known in dimension $g\geq 3$. 

In this paper we study various aspects of sign involutions, both of arithmetic and geometric nature.
Our first goal   is to investigate the \emph{existence of sign involutions} $\sigma$ on
twisted forms $X$ of   abelian varieties $A$, over general ground fields $k$
of arbitrary characteristic $p\geq 0$. These $\sigma$ are involutions on $X$ that become a sign involutions 
with respect to a suitable group law that arises on  some base-change. These varieties are usually introduced as torsors over
some abelian variety.
The following alternative point of view, developed  in  \cite{Laurent; Schroeer 2024} and \cite{Schroeer 2024},  is most suitable:
A \emph{para-abelian variety} is a proper scheme $X$ such that $X\otimes k'$ admits the structure of an abelian
variety, for some field extension $k\subset k'$.
It then   turns out that the  the subgroup scheme $A\subset\Aut_{X/k}$ that acts trivially on the 
numerically trivial part $\Pic^\tau_{X/k}$ is an abelian variety, 
and that the canonical $A$-action on $X$ is free and transitive. In turn, one may view the scheme $X$ as a torsor with respect to  the
abelian variety $A$ (the traditional point of view),
and obtains a   class $[X]$ in the \emph{Weil--Ch\^{a}telet group} $H^1(k,A)$.   
Our first main result relates these cohomology classes with the kernel $A[2]$ for the multiplication-by-two map
and the existence of sign involutions on $X$:

\begin{maintheorem}
(See Thm.\ \ref{divisibility relations})
Let $X$ be a para-abelian variety. Then the following are equivalent:
\begin{enumerate}
\item There is a sign involution $\sigma:X\ra X$.
\item We have $2\cdot [X]=0$ in the \emph{Weil--Ch\^{a}telet group} $H^1(k,A)$.
\item There is an  torsor $P$ with respect to $H=A[2]$ such that $X\simeq P\wedge^H A$.
\end{enumerate}
\end{maintheorem}

Here  $P\wedge^H A$ denotes the quotient of $P\times A$ by the diagonal $H$-action, usually  called
\emph{contracted product} or  \emph{associated fiber bundle}.
The main idea  idea for the above result is to introduce the \emph{scheme of sign involutions} $\Inv_{X/k}^\text{sgn}\subset\Aut_{X/k}$,
analyze the effect of the  conjugacy action on this subscheme, and derive consequences
using the general machinery of  twisted forms and non-abelian cohomology.

We now turn to   more geometric aspects:
Given an abelian variety $A$ with its standard sign involution $\sigma(x)=-x$, one can form the quotient 
$B=A/G$ with respect to the cyclic group  $G=\{e,\sigma\}$ of order two. This brings us into the realm of \emph{geometric invariant theory}:
Locally, the quotient arises from the ring of invariants in suitable coordinate rings for the abelian variety.
In characteristic two, not much seems to be known on the resulting proper normal scheme, and it would be highly interesting
to construct and understand a resolution of singularities.
Our second main result, which  is concerned with the numerically trivial part  $\Pic^\tau_{B/k}$  of the Picard scheme, which could shed some
  light on the problem:

\begin{maintheorem}
(See Thm.\ \ref{trivial pic-tau})
In the above situation,  the group scheme   $\Pic^\tau_{B/k}$ is trivial. 
\end{maintheorem}

This relies on    Grothendieck's two 
spectral sequences abutting to equivariant cohomology groups  \cite{Grothendieck 1957}. 
The   result is not difficult in the tame case $p\neq 2$, but requires a careful analysis  in the wild case $p=2$.  
Also note that  the statement immediately carries over to para-abelian varieties.
In dimension $g=1$ the para-abelian varieties $X$ are usually called \emph{genus-one curves}; we like to call them
\emph{para-elliptic curves}. These play an important role in the geometry and arithmetic of elliptic surfaces, in particular for bielliptic surfaces,
which also go by the name of hyperelliptic surfaces.
The above result shows  that the quotient by any sign involution is a \emph{Brauer--Severi curve}, that is,
a twisted form of $\PP^1$. 

Our third main result deals with the converse situation: Suppose there
is a degree-two morphism $f:X\ra B$ from a para-elliptic curve $X$ to some Brauer--Severi curve $B$. 
Then the projectivization $S=\PP(\shE)$
of the rank-two sheaf $\shE=f_*(\O_X)$ is a twisted form of a Hirzebruch surface with invariant $e=2$, and comes with 
a contraction  to a normal surface $S'$, having a unique singularity, which is often factorial.
The geometry of the situation is as follows:

\begin{maintheorem}
(See Section \ref{Morphisms})
Assumptions as above. Then   $f:X\ra B$ is the quotient by some sign involution $\sigma$ on the para-elliptic curve
$X$, and the latter  embeds into both surfaces $S$ and $S'$ as an anti-canonical curve. Moreover, $S'$ is the
anti-canonical model of $S$, and also a twisted form of the weighted projective space $\PP(1,1,2)$.
\end{maintheorem}

We also show that if there  are two different sign involutions $\sigma_1\neq \sigma_2$, the ensuing diagonal map gives an embedding
$X\subset B_1\times B_2$ into a product of Brauer--Severi curves. Such products where studied by 
Koll\'ar \cite{Kollar 2005} and Hogadi \cite{Hogadi 2009}. Again $X$ becomes an anti-canonical curve,
and it turns out that $B_1\times B_2$ embeds into $\PP^3$ if and  only if the factors are isomorphic.

\medskip
The paper is structured as follows:
In Section \ref{Scheme sign involutions} we recall the theory of para-abelian varieties $X$,
introduces the scheme of sign involutions $\Inv_{X/k}^\text{sgn}\subset\Aut_{X/k}$,
analyze the conjugacy action, and establish the link between sign involutions, cohomology classes, and 
structure reductions.
Section \ref{Picard scheme} is devoted to the Picard scheme  of the quotient $B=A/G$ of an abelian variety $A$
of arbitrary dimension $g\geq 0$ by a sign involution.
In Section \ref{Morphisms} we consider the case $g=1$, and unravel the geometry attached to 
degree-two maps $X\ra B$ from a para-elliptic curve $X$ to a Brauer--Severi curve $B$.

\begin{acknowledgement}
We like to thank the referees 
and Cec\'ilia Salgado
for their remarks, which helped to improve the paper.
The research was   conducted       in the framework of the   research training group
\emph{GRK 2240: Algebro-Geometric Methods in Algebra, Arithmetic and Topology}. The first two authors
where financially supported by the Deutsche For\-schungs\-gemeinschaft with a PhD grant in  GRK 2240/1, the first author
also with a PhD grant in  GRK 2240/2.
\end{acknowledgement}

%===========================================================
\section{The scheme of sign involutions}
\mylabel{Scheme sign involutions}

Let $k$ be a ground field of characteristic $p\geq0$, and $X$ be a proper scheme.
Then the group scheme $\Aut_{X/k}$ is locally of finite type, and the connected component
$\Aut_{X/k}^0$ of the neutral element $e=\id_X$ is of finite type (\cite{Matsumura; Oort 1967}, Theorem 3.7).  
By the Yoneda Lemma, the map $\sigma\mapsto\sigma^2$ defines a morphism of the scheme $\Aut_{X/k}$ to itself, which usually disrespects
the group law.
The  \emph{scheme of involutions} $\Inv_{X/k}$ is defined via a cartesian diagram
$$
\begin{CD}
\Inv_{X/k}	@>>>	\Aut_{X/k}\\
@VVV		@VV\sigma\mapsto\sigma^2 V\\
\Spec(k)	@>>e>	\Aut_{X/k}.
\end{CD}
$$
It contains the neutral element  and is stable under the
inverse map $\sigma\mapsto\sigma^{-1}$, but otherwise carries no further structure in general.

Now suppose that  $X$ can be endowed with the structure of an abelian variety. 
Recall that for each rational point $x_0\in X$, there is a unique group law
that turns $X$ into an abelian variety, with origin $0=x_0$.
Fix such a datum, and write $A$ for the abelian variety obtained by endowing  $X$ with the ensuing group law.
Note that $A$ can also be regarded as the pair $(X,x_0)$.
The automorphism group scheme becomes a semidirect product
$$
\Aut_{X/k}=A\rtimes\Aut_{A/k},
$$
where the normal subgroup on the left acts on $X$  by translations $x\mapsto a+x$.
The cokernel $\Aut_{A/k}$ on the right is an \'etale group scheme with countably many points, acting on   $A$ in the canonical way. 
Its  rational points are the automorphisms $\sigma:X\ra X$ fixing the origin $x_0$.
It  contains a canonical element, namely the  \emph{standard sign involution}  $x\mapsto -x$.
This defines a morphism $(-1):\Spec(k)\ra\Aut_{A/k}$. 
Its fiber with respect to  the canonical projection $A\rtimes\Aut_{A/k}\ra\Aut_{A/k}$
is denoted by  $A\otimes\kappa(-1)$.

\begin{lemma}
\mylabel{invariance}
The closed subscheme $A\otimes\kappa(-1)\subset \Aut_{X/k}$
is invariant under the conjugacy action of $\Aut_{X/k}$,   lies inside $\Inv_{X/k}$,
and does not  depend on the choice of the origin $x_0\in X$. 
\end{lemma}

\proof
Let $x,a,b\in A(R)$ and $\varphi\in\Aut_{A/k}(R)$ be   $R$-valued points, for some $k$-algebra $R$.
Then $x\mapsto a-x$ is some $R$-valued point of $A\otimes\kappa(-1)$. 
Conjugation by $(b,\id)$   is  
\begin{equation}
\label{conjugacy action}
x\longmapsto -b+x \longmapsto a-(-b+x)\longmapsto  (a + 2b)-x,
\end{equation}
whereas conjugation by  $(0,\varphi)$ takes the form 
$$
x\longmapsto \varphi^{-1}(x)\longmapsto a-\varphi^{-1}(x) \longmapsto  \varphi(a)-x.
$$
Both are $R$-valued points of $A\otimes\kappa(-1)$.
Furthermore, the   composition $x\mapsto a-x\mapsto a-(a-x)$ is the identity.
With the Yoneda Lemma, we see that $A\otimes\kappa(-1)$ is invariant under   conjugacy, and must be contained
in $\Inv_{X/k}$. 

Now let   $a_0\in X$ be another origin. 
The ensuing new group law and  negation are given by
$$
x\oplus y=x+y-x'_0\quadand \ominus x = -x + 2a_0,
$$ 
and thus $a\ominus x   =   (a +a_0)-x$. This shows that the closed subscheme $A\otimes\kappa(-1)\subset\Aut_{X/k}$
does not depend on the choice of origin.
\qed
 
\medskip
Recall that a proper scheme $X$ is called a \emph{para-abelian variety}
if there is a field extension $k\subset k'$ such that the base-change $X'=X\otimes k'$
admits the structure of an abelian variety.  This notation was introduced and studied by Laurent and the third author 
\cite{Laurent; Schroeer 2024}.
According to loc.\ cit., Proposition 5.2, the  closed subscheme $A\subset\Aut_{X/k}$ that acts trivial on  $\Pic^\tau_{X/k}$ is an abelian
variety, and the canonical $A$-action on $X$ is free and transitive. The resulting class  
$$
[X]\in H^1(k,A)
$$
in the Weil--Ch\^{a}telet group is called the \emph{cohomology class} of the para-abelian variety. Note that since $A$ is smooth, the
\'etale and fppf topology yield the same cohomology groups
(\cite{GB III}, Theorem 11.7). Consequently, the   class $[X]$ has some finite order; this number is
usually called   \emph{period} $\per(X)\geq 1$.

Conversely, if $H$ is any  commutative  group scheme, with a torsor $P$ and a homomorphism $H\ra A$,
we get a para-abelian variety $X=P\wedge^H X_0$. The latter denotes the quotient of $P\times X_0$
by the diagonal action $h\cdot (p,x) = (h\cdot p,h+x)$, and $X_0$ is the underlying scheme of the abelian variety $A$.
By construction, this $X$ is a twisted form of $X_0$.

Recall that the \emph{index} $\ind(X)\geq 1$ is the greatest common divisor
of the degrees $[\kappa(a):k]$ for the closed points $a\in X$. This is indeed the index
for the image of the degree map $\CH_0(X)\ra\ZZ$ on the Chow group of zero-cycles.
Note that in dimension one this can also be seen as the degree map  on the Picard group.
According to \cite{Lang; Tate 1958}, Proposition 5  the divisibility property $\per(X)\mid\ind(X)$ holds, and both numbers have the same prime factors.
 
As explained in \cite{Schroeer; Tziolas 2024}, Section 3, the group scheme $\Aut_{X/k}$ is a twisted form of $\Aut_{X_0/k}$ with respect to the
conjugacy action. In turn, the conjugacy-invariant closed subscheme $A\otimes\kappa(-1)\subset\Aut_{X_0/k}$ becomes
a closed subscheme $$\operatorname{Inv}^\text{sgn}_{X/k}\subset\Aut_{X/k},$$ which we call the \emph{scheme
of sign involutions}. Any automorphism $\sigma:X\ra X$ belonging to  $\operatorname{Inv}^\text{sgn}_{X/k}$
is called a \emph{sign involution}. 

\begin{theorem}
\mylabel{divisibility relations}
For each para-abelian variety  $X$  of dimension $g\geq 0$, the following three conditions are equivalent:
\begin{enumerate}
\item There is a sign involution $\sigma:X\ra X$.
\item We have $2\cdot [X]=0$ in the \emph{Weil--Ch\^{a}telet group} $H^1(k,A)$.
\item There is an  torsor $P$ with respect to $H=A[2]$ such that $X\simeq P\wedge^H A$.
\end{enumerate}
It these conditions hold we have the divisibility property $ \ind(X)\mid 4^g$.
\end{theorem}
 
\proof
We start with some general observations:
The first projection $$\Aut_{X_0/k}=A\rtimes\Aut_{A/k}\lra A$$ identifies the scheme of sign involutions 
$Z_0=\Inv^\text{sgn}_{X_0/k}=A\otimes\kappa(-1)$ with a copy of $X_0=A$. According to \eqref{conjugacy action}, the kernel for the
conjugacy homomorphism 
$A\ra\Aut_{Z_0/k}$   is $A[2]$, so this
factors over   multiplication-by-two map  $ A\stackrel{2}{\ra} A$.
It is now convenient to write  $X=T\wedge^AX_0$ for some $A$-torsor $T$.
Note that since the $X_0$ is the trivial $A$-torsor, one actually has $T=X$.
What is important now is that  the scheme of sign involutions
$Z=\Inv^\text{sgn}_{X/k} $ coincides with $Z=T\wedge^AZ_0$, and
the latter is the quotient of $T\times Z_0$ by  the $A$-action $a\cdot (t,z_0)=(a+t, 2a+z_0)$.

This quotient   can be computed as  successive quotients, first for the action of  $H=A[2]$ and
then for the induced action of $A/A[2]$. The group $H$ acts trivially on the second factor, hence
$H\backslash(T\times X_0)=(H\backslash T)\times X_0$. 
In light of the short exact sequence
\begin{equation}
\label{multiplication by two}
0\lra H\lra A\stackrel{2}{\lra} A\lra 0,
\end{equation}
we may regard $\bar{T}=H\backslash T$ as the $A$-torsor induced from $T$ with respect to $A\stackrel{2}{\ra} A$.
In other words  $Z=\bar{T}\wedge^{\bar{A}} Z_0$, where we write $\bar{A}=A/H=A$ to indicate the nature of the action.
By construction, the $\bar{A}$-action on $Z_0$ is free and transitive, so the projection $\bar{T}\otimes\kappa(-1)\ra Z$
is an isomorphism. 
We conclude that there is a rational point $\sigma\in Z$ if and only if the torsor $\bar{T}$ is trivial.
 
From the short exact sequence \eqref{multiplication by two} 
we get  a long exact sequence
$$
H^0(k,A)\stackrel{2}{\lra} H^0(k,A)\lra H^1(k,H)\lra H^1(k,A)\stackrel{2}{\lra} H^1(k,A).
$$
It follows that the element $[X]=[T]$ in $H^1(k,A)$ is annihilated by two if and only 
if there is an $H$-torsor $P$ such that such that $X\simeq P\wedge^HX_0$, giving the equivalence of (ii) and (iii).
Similarly, we see that $[X]=[T]$ is annihilated by two if and only if $\bar{T}$ is trivial.
Together with the previous paragraph this gives the equivalence of (i) and (ii).

It remains to verify the divisibility property of the index. This is just a special case of general fact:
Suppose $X$ has period $n\geq 1$. From the long exact sequence for the multiplication-by-$n$ map we see
that the quotient of $X$ by $A[n]$ contains a rational point, so its fiber $Z\subset X$ is a torsor with respect to $A[n]$.  According to \cite{Mumford 1970}, page 147 
the kernel  $A[n]$ is finite of length $l=n^{2g}$. Clearly, the torsor $Z$ has the same length, hence $X$ contains
a zero-cycle of degree $n^{2g}$. Now if (ii) holds, we have $n\mid 2$, and thus $\ind(X)\mid 4^g$.
\qed

\medskip
Recall that for each $m\geq 1$ there is an identification $H^1(k,\mu_m)=k^\times/k^{\times m}$.
Suppose now that $k$ contains a primitive $m$-th root of unity, such that $\mu_n\simeq(\ZZ/m\ZZ)_k$.
Let us recall the following result of Lang and Tate (\cite{Lang; Tate 1958}, Theorem 8):
Assume that  the ground field $k$, the abelian variety $A$, and the integer $m\geq 0$ satisfies the following conditions:
The   $\ZZ/m\ZZ$-module  $k^\times/k^{\times m}$ contains a free module of infinite rank, the quotient $A(k)/mA(k)$ is finite,
and  $A(k)$ contains an element of order $m$. Then  the Weil--Ch\^{a}telet group $H^1(k,A)$ contains infinitely many
elements $X$ whose period and index equals $m$.
Note that for global fields $k$, the first two conditions are automatic, and the third can be obtained after a finite extension,
provided the abelian variety has dimension $g\geq 1$ and the characteristic exponent $p\geq 1$ of $k$ is prime to $m$.

%===========================================================
\section{The Picard scheme of the quotient}
\mylabel{Picard scheme}

Let $A$ be an abelian variety, with its standard sign involution $\sigma(x)=-x$.
Write $G\subset\Aut(A)$ the corresponding subgroup of order two.
The quotient $B=A/G$ is a projective scheme  that is geometrically integral and geometrically normal, with $h^0(\O_B)=1$.
Following \cite{Fanelli; Schroeer 2020a}, Section 2, we write $\Sing(B/k)$ for the \emph{locus of non-smoothness}.
In contrast to the   \emph{locus of non-regularity} $\Sing(B)$, it comes with a scheme structure,
defined via Fitting ideals for K\"ahler differentials.

Let $\Pic^\tau_{B/k} $ be the open-and-closed subgroup scheme inside the Picard scheme comprising numerically trivial invertible sheaves.
Its Lie algebra is $H^1(B,\O_B)$, and the group scheme of connected components is the torsion part of the N\'eron--Severi group scheme.
It therefore encodes important information on $B$.  In dimension two,  $B=A/G$ yields the classical Kummer surfaces,
which give rise to  K3 surfaces, and in characteristic $p=2$ also to rational surfaces (\cite{Shioda 1974} and \cite{Katsura 1978}).
In both cases the tau-part of the Picard scheme vanishes. This generalizes to higher dimensions:

\begin{theorem}
\mylabel{trivial pic-tau}
The group scheme $\Pic^\tau_{B/k}$ is trivial. Moreover,  $\Sing(B/k)$ is finite, 
and is contained in the image of the fixed scheme $A^\sigma=A[2]$. 
\end{theorem}

\proof
It suffices to treat the case that $k$ is algebraically closed.
Write $q:A\ra B$ for the quotient map, let $U\subset A$ be the complement of the fixed scheme $A^\sigma=A[2]$, and $V=q(U)$ be its image.
The induced map $q:U\ra V$ is a $G$-torsor, in particular smooth.
According to \cite{EGA IVd}, Theorem 17.11.1 the smoothness of $U$ ensures the smoothness of $V$.
Thus $\Sing(B/k)$ is contained in the image of $A[2]$, and is therefore finite.

The structure sheaf $\O_A$ has a $G$-linearization, and thus comes with \emph{equivariant cohomology groups}
$H^i(A,G,\O_A)$, and likewise we have $H^i(A,G,\O_A^\times)$.
According to \cite{Grothendieck 1957}, Section 5.2, for every abelian sheaf $F$ on $A$ endowed with a $G$-linearization
there are  two spectral sequences
\begin{equation}
\label{grothendieck spectral sequence}
E_2^{rs} = H^r(G,H^s(A,F) )\quadand E_2^{rs} = H^r(B, \uH^s(G,F)),
\end{equation}
both with equivariant cohomology   $H^{r+s}(A,G,F)$ as abutment. For $F=\O_A^\times$ this gives two    exact sequences 
\begin{equation}
\label{five term sequences}
\begin{tikzcd} [row sep=tiny, column sep=small ]
0\ar[r]	& \Pic(B)\ar[r]	& H^1(A,G,\O_A^\times)\ar[r]\ar[d,equal]		& H^0(B,P)\ar[r]	& H^2(B,\O_B^\times)\\	
0\ar[r]	& H^1(G,k^\times)\ar[r]	& H^1(A,G,\O_A^\times)\ar[r]		& \Pic(A)^G\ar[r]		& H^2(G,k^\times),
\end{tikzcd}
\end{equation}
% \begin{equation}
% \label{five term sequences}
% \begin{tikzcd} [row sep=tiny, column sep=small ]
% 0\ar[r]	& \Pic(B)\ar[dr]	&				& \Pic(A)^G\ar[r]	& H^2(G,k^\times)\\
% 	&		& H^1(A,G,\O_A^\times)\ar[ur]\ar[dr]\\	
% 0\ar[r]	& H^1(G,k^\times)\ar[ur]			&			& H^0(B,F)\ar[r]	& H^2(B,\O_B^\times),
% \end{tikzcd}
% \end{equation}
where the abelian sheaf $P=\uH^1(G,\O_A^\times)$  is supported by the singular locus of $B$, and the composition
$\Pic(B)\ra H^1(A,G,\O_A^\times)\ra \Pic(A)^G$ is given by pullback of invertible sheaves.
Recall that the cohomology groups for the cyclic group $G=\{e,\sigma\}$ are given by
$$
H^{2j+1}(G,M) = \frac{\Kernel(\sigma+\id)}{\Image(\sigma-\id)} \quadand H^{2j+2}(G,M)= \frac{\Kernel(\sigma-\id)}{\Image(\sigma+\id)},
$$
for any $G$-module $M$.
It follows that $H^{2}(G,k^\times)$ vanishes,  because $G$ acts trivially on $k^\times$,  and $k^\times=k^{\times2}$, whereas
$H^1(G,k^\times)=\mu_2(k)=\{\pm 1\}$.
According to \eqref{five term sequences}, the kernel for  $\Pic(B)\ra\Pic(A)$  is the intersection 
of $\Pic(B)\cap H^1(G,k^\times)$ inside the equivariant cohomology group. Furthermore, the image of $\Pic^\tau(B)\ra\Pic(A)$   is contained in 
$$
\Pic^\tau(A)\cap\Pic(A)^G=A(k)\cap \Pic(A)^G =A(k)[2]=\Pic(A)[2].
$$
This already shows that the group  scheme $\Pic^\tau_{B/k}$ must be finite. 
It also settles the case of dimension $g=1$:  Then $B$ is a normal curve with
finite Picard scheme. The latter is smooth, according to \cite{Mumford 1966}, Section 27 because $H^2(B,\O_B)=0$.
Consequently $B=\PP^1$, and thus $\Pic^\tau_{B/k}=0$. 

From now on, we assume that we are in dimension $g\geq 2$.  At each $a\in A[2]$, the induced $G$-action on the local ring $\O_{A,a}$
is ramified only at the origin. It follows that the local ring at the image $b\in B$ is singular,   and that the finite degree-two extension
$\O_{B,b}\subset\O_{A,a}$   is not flat: the arguments  in  \cite{Lorenzini; Schroeer 2020}, last paragraph in the proof for Proposition 3.2,
hold true for the action of our group $G$ of order two in characteristic $p\geq 0$.
Consequently, the quotient map $q:A\ra B $ induces
a bijection between $A[2]$ and $\Sing(B)$.
Furthermore, the short exact sequence $0\ra \O_B\ra q_*(\O_A)\ra \shF\ra 0$ defines a coherent sheaf  $\shF$ 
that is invertible on the open set $V=\Reg(B)$, but not at the points $b\in\Sing(B)$. 

We claim that the canonical map $\Pic(B)\ra \Pic(A)^G$ is injective. Equivalently,
the   intersection $\Pic(B)\cap H^1(G,k^\times)$ inside $H^1(A,G,\O_A^\times)$ is trivial. The group $H^1(G,k^\times)=\mu_2(k)$ vanishes
in characteristic two, so only the  case $p\neq 2$ requires attention.
Then the \emph{trace map} $q_*(\O_A)\ra\O_B$, which sends a local section viewed as an $\O_B$-linear homothety  to its trace,
 gives a splitting $q_*(\O_A)=\O_B\oplus\shF$, thus $\shF$ satisfies Serre's Condition  $(S_2)$.
The canonical identification $ \shF_V\otimes\shF_V^\vee=\O_V$ yields an element in 
$\Gamma(V,q_*(\O_A)\otimes\shF^\vee)=\Gamma(U,q^*(\shF^\vee))$ without zeros, and it 
follows that the invertible sheaf  $\shF|V$   becomes   trivial   on $U$.
Using the diagram \eqref{five term sequences} for the  quotient $V=U/G$ instead of $B=A/G$, we conclude that $\shF|V$
generates the kernel of $\Pic(V)\ra\Pic(U)$. Seeking a contradiction, we now assume that there is a non-trivial invertible sheaf $\shL$ on $B$
that becomes trivial on $A$, we therefore must have $\shL|V=\shF|V$.
Using that both $\shL$ and $\shF$ satisfies Serre's Condition $(S_2)$ together with    \cite{Hartshorne 1994}, Theorem 1.12
we infer that $\shL=\shF$, contradicting that $\shF$ is not invertible. This establishes our claim.
In turn, the canoncial map $\Pic^\tau(B)\ra\Pic(A)[2]$ becomes an inclusion $\Pic^\tau(B)\subset\Pic(A)[2]$.

We next check that for $p\neq 2$ the finite group scheme $\Pic^\tau_{B/k}$ is reduced. Equivalently, its Lie algebra $H^1(B,\O_B)$ vanishes.
To see this,  consider the spectral sequences \eqref{grothendieck spectral sequence}
with the additive sheaf $\O_A$ instead the multiplicative sheaf  $\O_A^\times$.
For $i\geq 1$, the vector spaces  $H^i(G,k)$ are annihilated by the group order $|G|=2$. For $p\neq 2$ they consequently    vanish, and we obtain inclusions
$$
H^1(B,\O_B)\subset H^1(A,G,\O_A)\subset H^1(A,\O_A)^G.
$$
Moreover, the term on the right also vanishes because $G$ acts via the sign involution on the cohomology group, according
(\cite{Roessler; Schroeer 2022}, proof of Proposition 2.3). This establishes the claim.

To proceed we use the fact that for any finite commutative group scheme $N$ the
isomorphism classes of $N$-torsors $B'\ra B$ corresponds to homomorphisms of group schemes  $N^*\ra \Pic_{B/k}$,
where $N^*=\uHom(N,\GG_m)$ denotes the \emph{Cartier dual} (see \cite{Raynaud 1970}, Proposition 6.2.1,
and also the discussion in  \cite{Schroeer 2017}, Section 4).

The constant group scheme $N=(\ZZ/2\ZZ)_k$ has Cartier dual $N^*=\mu_2$.
Suppose we have an inclusion $\mu_2\subset\Pic^\tau_{B/k}$ such that the composite
map $\mu_2\ra\Pic^\tau_{A/k}$ remains a monomorphism.
The corresponding   $N$-torsor $B'\ra B$ thus induces a non-trivial $N$-torsor $A'\ra A$.
According to the Serre--Lang Theorem (\cite{Mumford 1970}, page 167), there is a unique structure of an abelian variety for $A'$ so that
$A'\ra A$ is a homomorphisms. This gives an embedding $N\subset A'$ defined by  a 2-division point $a'\in A'$.
The composite $A' \ra B$ is the quotient for the action of  $N\rtimes\{\pm1\}$.
 Since this semidirect product is actually   a direct product, the projection $A'\ra B'$ must be the quotient by
$G=\{\pm1\}$. Now choose a closed point $x'\in A'$ with $2x'=a'$. It follows that
the orbit $G\cdot x'=\{\pm x'\}$, viewed as a rational point on $B'$, is fixed by the the $N$-action,
contradiction. This settles the case $p\neq 2$:
Then $\mu_2=(\ZZ/2\ZZ)_k$, and  we see that $\Pic^\tau(B)\subset\Pic(A)[2]$ is trivial.
We already saw in the previous paragraph that $\Pic^\tau_{B/k}$ is reduced, and infer that it must be trivial.
 
It remains to treat the case $p=2$, where the arguments in some sense  run  parallel to the preceding paragraph. At each $a\in A[2]$,
the local ring at the image $b\in B$ is singular, with   $\depth(\O_{B,b})=2$, according to \cite{Lorenzini; Schroeer 2020}, Proposition 3.2.
Note that this is in stark contrast to the situation $p\neq 2$, when such rings of invariants are Cohen--Macaulay. 
Again we consider the short exact sequence
$0\ra\O_B\ra q_*(\O_A)\ra\shF\ra 0$ of coherent sheaves on $B$.   
%According to \cite{Ekedahl 1988}, Section 1  we have $\shF|V=\O_V$.
For the images $b\in B$ of the $a\in A[2]$, the short exact sequence of local cohomology
$$
H^0_b(B,q_*(\O_A))\lra H^0_b(B,\shF)\lra H^1_b(B,\O_B),
$$
reveals that $H^0_b(B,\shF)=0$, in other words, $\shF$ is torsion-free. 
The trace map $q_*(\O_A)\ra\O_B$ vanishes on the subsheaf $\O_B\subset q_*(\O_A)$ since we are in characteristic two.
The induced map  $\shF\ra\O_B$ is bijective on the locus where $\shF$ is invertible, which
one easily sees by a local computation.
This gives an inclusion $\shF\subset\O_B$. 
%So the adjunction map $\shF\ra i_*(\shF|V)=\O_B$ is injective, hence $\shF$ is   a sheaf of ideals.
Using that $\shF$ is not invertible we infer $H^0(B,\shF)=0$.
The exact sequence
$$
H^0(B,\shF)\lra H^1(B,\O_B)\lra H^1(A,\O_A) 
$$
ensures that the map on the right  is injective. On the other hand, its kernel
is the Lie algebra for the kernel of $\Pic^\tau_{B/k}\ra \Pic_{A/k}[2]$.
It follows that this map is actually a closed embedding $\Pic^\tau_{B/k}\subset\Pic_{A/k}[2]$.

Now we use that the Lie algebra of any group scheme in characteristic $p>0$
carries as additional structure the  $p$-map $x\mapsto x^\pmap$ and becomes
a \emph{restricted Lie algebra} (see \cite{Schroeer; Tziolas 2024}, Section 1 for more details).
Suppose $H^1(B,\O_B)\neq 0$. Then there is a $p$-closed vector $x\neq 0$, in other words $x^{[p]}$ is a multiple of $x$.
The case $x^\pmap\neq 0$ yields an inclusion of $\mu_p\subset B$ where the composite map $\mu_p\ra A$ is injective.
We saw above that this is impossible.
In turn we must have $x^\pmap=0$. This gives an inclusion of $N^*=\alpha_p$ into $B$ where the composite map 
$\alpha_p\ra A$ remains injective. The Cartier dual is $N=\alpha_p$.
Thus we get a non-trivial $\alpha_p$-torsor $B'\ra B$ for $\alpha_p$ whose base-change $A'\ra A$ remains non-trivial.
A similar situation with $ N^*=(\ZZ/2\ZZ)_k$ and    $N=\mu_p$ arise if there is a point of order two on $\Pic_{B/k}$.
In both cases the discussion in \cite{Roessler; Schroeer 2022}, beginning of Section 2 shows that
 $A'$ has the structure of an abelian variety so that the projection $A'\ra A$ is a homomorphism, and we get an inclusion $N\subset A'$.
The composition $A'\ra B$ is the quotient by the group scheme $N\rtimes\{\pm 1\}$. Again this  is actually a direct product.
In the cartesian diagram
$$
\begin{CD}
A'	@>>>	B'\\
@VVV		@VVV\\
A	@>>>	B
\end{CD}
$$
the vertical maps are quotients by the action of the infinitesimal group scheme $N$, and the horizontal maps are quotients
by $G=\{\pm1\}$. Fix some $a'\in A'[2]$, with image $b'\in \Sing(B')$, and consider the ring of invariants
$\O_{B',b'}\subset\O_{A',a'}$.
According to \cite{Lorenzini; Schroeer 2020}, Lemma 3.3 no element $f\in\maxid_{ a'}\smallsetminus\maxid_{a'}^2$ is $G$-invariant.
It follows that the infinitesimal neighborhood $\Spec(\O_{a'}/\maxid_{a'}^2)$ maps to
$Z'=\Spec(\O_{b'}/\maxid_{b'})$, and therefore  the same holds for the  orbit $N\cdot\{a'\}$.
In light of the above commutative diagram, the $N$-action on $B'$
is not free, contradiction. 
\qed

\medskip
The result immediately carries over to para-abelian varieties, because the formation of both the quotient $B=A/G$ and the Picard scheme $\Pic_{B/k}$
commutes with ground field extensions.
The para-abelian varieties $X$ of dimension $g=1$ are usually called \emph{genus-one curves}.
Throughout, we shall prefer the term \emph{para-elliptic curves}. 
These are   twisted forms of   elliptic curves. The moduli stack of such curves was studied
by the second author \cite{Dang 2022}.
Recall that the \emph{Brauer--Severi varieties} $Y$ are twisted forms
of projective  space $\PP^n$, for some $n\geq 0$. For more details we refer to \cite{Artin 1982}.
In case  $n=1$ we also say that $Y$ is a  \emph{Brauer--Severi curve}. 

\begin{corollary}
\mylabel{brauer--severi quotient}
Assumption as in the proposition, and suppose additionally $g=1$. Then the corresponding quotient $B=X/G$ is a Brauer--Severi curve.
\end{corollary}

\proof
The scheme $B$ is geometrically normal and of dimension one, hence smooth.
According to the theorem, the Picard scheme is discrete.
It follows that the tangent space $H^1(B,\O_B)$   vanishes.
If there is a rational point $a\in X$, the resulting invertible sheaf $\shL=\O_B(a)$ is very ample,
with $h^0(\shL)=2$, and we obtain an isomorphism $B\ra\PP^1$. 
\qed

\medskip
In dimension $g=2$ and characteristic $p\neq 2$, the quotient $B=A/\{\pm 1\}$ is called a \emph{Kummer surface},
and is a K3 surface with rational double points.
For $p=2$, the quotient $B$ is either a K3 surface with rational double points,
or a rational surface with an elliptic singularity. This was discovered by Shioda \cite{Shioda 1974},
see also \cite{Katsura 1978},  \cite{Schroeer 2007}, \cite{Schroeer 2009} and \cite{Kondo; Schroeer 2019}.
The formation of such quotients is studied by the first author \cite{Bergqvist 2023}.
Little seems to be know on the quotient in higher dimensions, in particular in characteristic two, 
compare Schilson's  investigation \cite{Schilson 2018}, \cite{Schilson 2020}.

%===========================================================
\section{Morphisms to Brauer--Severi curves}
\mylabel{Morphisms}

Let $X$ be a para-elliptic curve over a ground field $k$. If there is a sign involution $\sigma:X\ra X$,
the quotient $B$ by the corresponding group of order two
is a Brauer--Severi curve, according  to Corollary \ref{brauer--severi quotient}. 
In this section we conversely assume that our para-elliptic curve $X$ admits 
a morphism $f:X\ra B$ of degree two to some Brauer--Severi curve $B$, and derive several geometric consequences.

First note that the corresponding function field extension $k(B)\subset k(X)$ has degree two. It must be separable, because 
 $X$ and $B$ are smooth of different genus. So this is a Galois extension, and the Galois group $G$
is cyclic of order two. Let $\sigma\in G$ be the generator.

\begin{proposition}
\mylabel{automorphism sign involution}
The automorphism $\sigma:X\ra X$ is a sign involution.
\end{proposition}

\proof
It suffices to treat the case that $k$ is algebraically closed. The action is not free, because
$\chi(\O_X)=0\neq 2=|G|\cdot\chi(\O_B)$. Choose a fixed point $x_0\in X$, and regard  $E=(X,x_0)$ as an elliptic curve.
If   $\Aut(E)$ is   cyclic, there is a unique element of order two, and we infer that $\sigma$ equals the sign involution.
Suppose now that $\Aut(E)$ is non-cyclic. According to \cite{Deligne 1975}, Proposition 5.9
this  group    is either the semi-direct product $\ZZ/3\ZZ\rtimes\mu_4(k)$ in characteristic $p=3$,
or $Q\rtimes\mu_3(k)$ in characteristic $p=2$, where $Q=\{\pm 1,\pm i,\pm j, \pm k\}$ denotes the quaternion group.
In these groups, the respective elements
 $(0,-1)$ and $(-1, 1)$ are the only ones  of order two, and we again conclude that $\sigma $ coincides with the sign involution.
\qed

\begin{proposition}
\mylabel{cokernel sheaf}
The cokernel for the inclusion $\O_B\subset f_*(\O_X)$ is isomorphic to $\omega_B$,
and the resulting extension $0\ra\O_B\ra f_*(\O_X)\ra\omega_B\ra 0$  of coherent sheaves splits.
\end{proposition}

\proof
The sheaf $f_*(\O_X)$ has rank two and is torsion-free, hence is locally free. The inclusion of $\O_B$ is locally
a direct summand, so the cokernel $\shL$ is invertible.
We have  $0=\chi(\O_X)=\chi(\O_B)+\chi(\shL)=2+\deg(\shL)$
and conclude $\deg(\shL)=-2$. Since $\deg:\Pic(B)\ra \ZZ$ is injective, this gives $\shL\simeq\omega_B$.
The extension yields a class in $\Ext^1(\omega_B,\O_B)=H^1(X,\omega_B^{\otimes-1})$, which vanishes by Serre Duality.
So the extension splits.
\qed

\medskip
Choose a splitting and set  $\shE=f_*(\O_X)=\O_B\oplus\omega_B$. 
The smooth surface 
$$
S=\PP(\shE)=\Proj(\Sym^\bullet\shE)
$$
 is a twisted form of the Hirzebruch surface
$S_0=\PP(\shE_0)$, where $\shE_0=\O_{\PP^1}\oplus\O_{\PP^1}(-2)$.
Let us call  $S$ the \emph{twisted Hirzebruch surface} attached to the
Brauer--Severi curve $B$.
Since $f:X\ra B$ is affine, the invertible sheaf $\O_X$ is relatively very ample,
and we get a closed embedding $X\subset S$.
By abuse of notation we also write  $f:S\ra B$ for the extension of   our original morphism on $X$.

Recall that each invertible quotient $\shE\ra \shN$
defines a section $s:B\ra S$, whose image $D$ has self-intersection $D^2=\deg(\shN)-\deg(\shN')$, where
$\shN'\subset\shE$ is the kernel.
For more details we refer to \cite{Fanelli; Schroeer 2020b}, Section 6.
In particular,   $\pr_1:\shE\ra\O_B$ yields
a curve $D\subset S$ with $D^2=2$,
whereas   $\pr_2:\shE\ra\omega_B$
gives some $E\subset S$ with $E^2=-2$, and the two sections are disjoint.
The Adjunction Formula gives $(\omega_S\cdot D)=-4$ and $(\omega_S\cdot E)=0$. Hence  $\omega_S=f^*(\omega_B^{\otimes 2})\otimes\O_S(-2E)$, because
both sides have the same intersection numbers with $D$ and $E$.
In particular
$c_1^2=(\omega_S\cdot\omega_S)= -8\cdot \deg(\omega_B) + 4\cdot E^2=  8$.
Setting 
$$
\omega_S^{\otimes 1/2}=f^*(\omega_B)\otimes\O_S(-E),
$$ 
we get an invertible sheaf whose square is isomorphic
to the dualizing sheaf. In other words, the surface $S$ comes with a canonical \emph{theta characteristic},
or \emph{spin structure}, compare \cite{Atiyah 1971} and \cite{Mumford 1971}.

\begin{proposition}
\mylabel{spin structure}
The dual sheaf $\shL=\omega_S^{\otimes -1/2}$ is globally generated with $h^0(\shL)=4$.
The image of the resulting   $r:S\ra\PP^3$ is an integral normal surface $S'\subset\PP^3$ of degree two,
and  the induced morphism $r:S\ra S'$ is the contraction of $E$. Moreover, the image $a=r(E)$ is a rational point,
the local ring $\O_{S',a}$  is singular, and the   restriction $r|X$ is a closed embedding.
\end{proposition}

\proof
Our sheaf  has 
intersection numbers $(\shL\cdot\shL)=2$ and $(\shL\cdot E)=0$. 
Serre Duality gives $h^2(\shL) = h^0(\omega_S^{\otimes 3/2})=0$, and Riemann--Roch yields
$$
h^0(\shL)\geq \chi(\shL) = \frac{c_1^2/4 + c_1^2/2}{2} +\chi(\O_S) =  (2+4)/2 + 1 = 4.
$$
The base locus $\Bs(\shL)$ is   contained in $E$, because $\omega_B^{\otimes-1}$ is globally generated.
The short exact sequence $0\ra f^*(\omega_B^{\otimes-1})\ra\shL\ra \shL|E\ra 0$ yields an exact sequence
$$
0\lra H^0(S,f^*(\omega_B^{\otimes-1})) \lra H^0(S,\shL)\lra H^0(E,\O_E),
$$
consequently $h^0(\shL)\leq h^0(\omega^{\otimes-1}_B)+h^0(\O_E)=4$. This ensures $h^0(\shL)=4$, and that  $\shL$ is globally generated.

In turn, our spin structure yields a morphism 
$r:S\ra \PP^3$
with $r^*(\O_{\PP^3}(1))=  \omega_S^{\otimes -1/2}$. It therefore contracts $E$. Moreover,
the image    $S'\subset\PP^3$ is integral and two-dimensional, of some degree $n\geq 1$.
 This image is not a  plane, because  the morphism is defined by the
complete linear system $H^0(S,\shL)$. From $2=(\shL\cdot\shL)=\deg(S/S')\cdot n$ we infer that $S\ra S'$ is birational 
and $n=2$. The Adjunction Formula gives $\omega_{S'}=\O_{S'}(-2)$, consequently $r^*(\omega_{S'})=\omega_S$.
It follows that the birational morphism $r:S\ra S'$ is in Stein factorization. Since $\Pic(S)$ has rank two,
the exceptional divisor is irreducible, whence must coincide with $E$.

The image $a=r(E)$ is a rational point, because $h^0(\O_E)=1$. The local ring $\O_{S',a}$ must be singular,
because otherwise $S=\Bl_a(S')$, such that $E=r^{-1}(a)$ must be a projective line with $E^2=-1$, contradiction.

It remains to verify that the curves $X,E\subset S$ are disjoint. Since $\deg(X/B)=2$ we have $\omega_S = \O_S(-X)\otimes f^*(\shN)$
for some invertible sheaf $\shN$ on $B$. The Adjunction Formula gives
$$
0=(\omega_S\cdot X)+ X^2 = -X^2 + 2\deg(\shN) + X^2.
$$
Consequently $\shN$ is trivial, and $\omega_S=\O_S(-X)$. This gives $X^2=c_1^2=8$, and furthermore $(X\cdot E) = -(\omega_S\cdot E)=0$.
Thus the integral curves $X$ and $E$ must be disjoint, hence $r|X$ is a closed embedding.
\qed

\medskip
Note that the local ring $\O_{S',a}$ is factorial provided that $B\not\simeq\PP^1$.
The above also shows that the image  $S'=r(S)$ can also be viewed as the \emph{anti-canonical model} $P(S,-K_S)$
of the scheme $S$, which is defined as the homogeneous spectrum of the \emph{anti-canonical ring}
$R(S,-K_S) = \bigoplus_{t\geq 0}H^0(S,\omega_S^{\otimes t})$.

Recall that the \emph{weighted projective space} $\PP(d_0,\ldots,d_n)$ is the homogeneous spectrum 
 of $k[U_0,\ldots,U_n]$, where the generators have degrees  $d_i=\deg(U_i)$.
 The case $d_0=\ldots=d_n=1$
gives back the standard projective space $\PP^n$.
Let us  say that a closed subscheme of a Gorenstein surface is an \emph{anti-canonical curve} if its sheaf of ideals is
isomorphic to the dualizing sheaf.

\begin{proposition}
\mylabel{anti-canonical model}
The anti-canonical model $S'=P(S,-K_S)$ is a twisted form of the weighted projective space $\PP(1,1,2)$.
Moreover, $X\subset S$ and the resulting inclusion $X\subset S'$ are anti-canonical curves.
\end{proposition}

\proof
It suffices to treat the case that $k$ is algebraically closed. 
We claim that $S'$ is defined inside $\PP^3=\Proj k[T_0,\ldots,T_3]$
by the equation $T_0^2-T_1T_2=0$, for a suitable choice of homogeneous coordinates.
The main challenge is the case $p=2$:  
According to \cite{Arf 1941}, Satz 2 our quadric $X\subset\PP^3$ must be defined by  an equation of the form
$$
\sum_{i=1}^r(\alpha_iX_i^2 + X_iY_i + \gamma_iY_i^2) + \sum_{j=1}^s \delta_jZ_j^2=0,
$$
with $1\leq 2r+s\leq 4$, and non-zero coefficients $\delta_j$. Since $k$ is algebraically closed,
we can make a change of variables and  achieve $\delta_j=1$, and  furthermore $\alpha_i=\gamma_i=0$.
For $s\geq 1$ the coordinate change $Z_1=Z_1'+\ldots+Z'_s$ reduces us to the case $s=1$.
One now immediately sees that only  for $r=s=1$ the quadric $S'\subset\PP^3$ is normal and singular, and setting
$T_0=Z_1$ and $T_1=X_1$ and $T_2=Y_1$ gives the claim.
For $p\neq 2$ our quadric can be defined by an equation of the form $\sum_{j=0}^3\delta_jZ_j^2=0$,
and one argues similarly.

Consider the graded ring $A=k[U_0,U_1,U_2]$ with weights $(1,1,2)$.
The Veronese subring $A^{(2)}$ is generated by the homogeneous elements $U_0U_1, U_0^2,U_1^2,  U_2$,
which satisfy the relation $(U_0U_1)^2=U_0^2\cdot U_1^2$. This gives a surjection
$$
k[T_0,T_1,T_2,T_3]/(T_0^2-T_1T_2)\lra A^{(2)},
$$
defined by the assignments
$T_0\mapsto U_0U_1$ and $T_1\mapsto U_0^2$ and $T_2\mapsto U_1^2$ and $T_3\mapsto U_2$. 
Both rings   are integral of dimension three. Using Krull's Principal Ideal Theorem, we infer that
the above surjection is bijective. The homogeneous spectrum of $A^{(2)}$ coincides with $\PP(1,1,2)=\Proj(A)$,
and by the above also with $S'$.

We already saw in the previous proof that $\omega_S=\O_S(-X)$, hence $X\subset S$ is an anti-canonical curve.
From the Theorem of Formal functions one infers $f_*(\omega_S)$ is invertible, and this ensures
that the direct image coincides with $\omega_{S'}$. Using $X\cap E=\varnothing$ we infer $\omega_{S'}=\O_{S'}(-X)$.
\qed

\medskip
Now suppose that we have \emph{two} morphism $B_1\stackrel{f_1}{\leftarrow} X\stackrel{f_2}{\ra} B_2$
to Brauer--Severi curves, with $\deg(X/B_i)=2$. According to Proposition \ref{automorphism sign involution}, they come 
from  sign involutions $\sigma_1$ and $\sigma_2$, respectively.  

\begin{proposition}
\mylabel{diagonal embedding}
If $\sigma_1\neq \sigma_2$, the diagonal morphism $i:X\ra B_1\times B_2$ is a closed embedding, and its
image is an anti-canonical curve.
\end{proposition}

\proof
Let $A\subset\Aut_{X/k}$ be the subgroup scheme that fixes $\Pic^\tau_{X/k}$.
As discussed in Section \ref{Scheme sign involutions}, this is an elliptic curve, and the action on the para-elliptic curve $X$ is free and transitive.
Moreover, the dual abelian variety is identified with  $\Pic^0_{X/k}$. But note that the principal polarization stemming from the origin also
gives $A=\Pic^0_{X/k}$.
We saw in the proof of Proposition \ref{invariance} that the  two rational points $\sigma_1,\sigma_2\in \Inv_{X/k}^\text{sgn}$ differ
by the action of some non-zero $a\in A(k)$. In other words,
$\sigma_2(x) = a+\sigma_1(x)$.
It follows that there is no rational point $x\in X$ with $\sigma_1(x)=\sigma_2(x)$.
In particular, the fixed schemes $X^{\sigma_1}$ and $X^{\sigma_2}$ are disjoint.

To proceed, we assume that $k$ is algebraically closed. Let $x\in X$ be a closed point and write $y=i(x)=(b_1,b_2)$.
The  inverse image $i^{-1}(y)$ is the   intersection  
of the fibers $f_1^{-1}(b_1)\cap f_2^{-1}(b_2)$. This is just the spectrum of $\kappa(x)$, by the previous paragraph.
According to \cite{EGA IVd}, Corollary 18.12.6 the finite morphism $i:X\ra B_1\times B_2$ is a closed embedding.

By construction, we have $\deg(X/B_1)=\deg(X/B_2)=2$.   Set $V=B_1\times B_2$. 
Its Picard scheme $\Pic_{V/k}$ can seen as the Galois module $\Pic(V\otimes k^\sep)=\ZZ\times\ZZ$,
compare the discussion in \cite{Schroeer 2023}, Section 1. Obviously, the elements $(2,0)$ and $(0,2)$ are fixed
by $\Gal(k^\sep/k)$, hence the whole Galois action is trivial, and 
thus  $\Pic_{V/k}=(\ZZ\times\ZZ)_k$ is a constant group scheme. 
The dualizing sheaf  $\omega_V=\pr_1^*(\omega_{B_1})\otimes \pr_2^*(\omega_{B_2})$
has class $(2,2)$,  and we infer $\omega_V=\O_S(-X)$.
\qed

\medskip
Note that $\omega_V$ is anti-ample, so  the smooth surface $V=B_1\times B_2$ coincides with its anti-canonical model $P(V,-K_V)$.
Products of Brauer--Severi curves were studied by Koll\'ar \cite{Kollar 2005} and Hogadi \cite{Hogadi 2009}. Let us close this paper with
the following   observation:

\begin{proposition}
\mylabel{embedding into product}
The surface $V=B_1\times B_2$ admits an embedding into $\PP^3$ if and only if $B_1\simeq B_2$.
\end{proposition}

\proof
The Picard scheme is given by $\Pic_{V/k}=(\ZZ\times\ZZ)_k$. The classes $(-2,0)$ and $(0,-2)$ come from the preimages of the invertible sheaves 
on $B_1$ and $B_2$, and thus belong to the subgroup $\Pic(V)\subset\Pic_{V/k}(k)$.

Suppose we have $V\subset\PP^3$, and write $d\geq 1$ for its degree. From $\omega_V=\O_V(d-4)$  we
get $8=(\omega_V\cdot\omega_V) = d(d-4)^2$, and thus $d=2$. In particular, $V$ admits the spin structure $\omega_V^{\otimes 1/2} =\O_V(-1)$.
The dual sheaf $\shL=\O_V(1)$ has $h^0(\shL)=4$, which easily follows from the short exact sequence
$0\ra\O_{\PP^3}(-1)\ra \O_{\PP^3}(1)\ra \shL\ra 0$. Choose some non-zero global section $s\neq 0$ from $\shL$,
and let $D\subset V$ the resulting effective Cartier divisor. Suppose $D$ is reducible. Since $\deg(D)=2$ we see that there are two
components. Since $\shL$ has class $(1,1)$ in $\Pic_{V/k}(k)$, it follows that $D=D_1+D_2$, where the summands are preimages
of rational points on $B_1$ and $B_2$, respectively.
Thus both Brauer--Severi curves are copies of $\PP^1$. 
Suppose now that $D$ is irreducible. Then $\deg(D/B_i)=1$, so the morphism $D\ra B_i$ are birational. By Zariski's Main Theorem,
it  must be an isomorphism, and therefore $B_1\simeq B_2$.

Conversely, suppose there is an isomorphism $h:B_1\ra B_2$. Its graph  
defines an effective  Cartier divisor $D\subset B_1\times B_2$ with class $(1,1)\in\Pic_{V/k}(k)$. Set $\shL=\O_V(D)$.
Passing to the algebraic closure of $k$, we get $\shL=\pr_1^*(\O_{\PP^1}(1))\otimes \pr_2^*(\O_{\PP^1}(1))$,
and compute $h^0(\shL)=4$. Moreover, $\shL$ is very ample, and thus    defines a closed embedding $X\subset \PP^3$.
\qed

\medskip
Given a sign involution $\sigma: X\ra X$ and a non-zero rational point $a\in A(k)$, we get another sign involution
$x\mapsto a+\sigma(x)$.
We see that the   situation $B_1\stackrel{f_1}{\leftarrow} X\stackrel{f_2}{\ra} B_2$ with $\sigma_1\neq \sigma_2$ appears if and only if
the set  $\Inv^\text{sgn}_{X/k}(k)$ is non-empty
and the group $A(k)$ is non-trivial.

%===========================================================


\begin{thebibliography}{ccccc}

\bibitem{Arf 1941}
C.\ Arf:
Untersuchungen \"uber quadratische Formen in K\"orpern der Charakteristik 2. I. 
J.\ Reine Angew.\ Math.\ 183 (1941), 148--167. 

\bibitem{Artin 1969}
M.\ Artin: 
Algebraization of formal moduli I. 
In:  
D.\ Spencer, S.\ Iyanaga (eds.),
Global Analysis, pp. 21--71.
Univ.\ Tokyo Press, Tokyo, 1969. 

\bibitem{Artin 1982}
M.\ Artin:
Brauer--Severi varieties.
In:  F.\ van Oystaeyen, A.\ Verschoren (eds.),
Brauer groups in ring theory and algebraic geometry, pp.\ 194–-210,
Springer, Berlin-New York, 1982. 

\bibitem{Atiyah 1971}
M.\ Atiyah:
Riemann surfaces and spin structures.
Ann.\ Sci.\ \'Ecole Norm.\ Sup.\  4 (1971), 47--62. 

\bibitem{Bergqvist 2023}
J.\ Bergqvist:
The Kummer constructions in families.
Dissertation, D\"usseldorf (2023)
\url{https://docserv.uni-duesseldorf.de/servlets/DocumentServlet?id=63084}.

\bibitem{Dang 2022}
T.\ Dang:
Cohomology of certain Artin stacks.
Dissertation, D\"usseldorf (2022),  
\url{https://nbn-resolving.org/urn/resolver.pl?urn=urn:nbn:de:hbz:061-20220822-084645-3}.

\bibitem{Deligne 1975}
P.\ Deligne:
Courbes elliptiques: formulaire d'apr\`es J.\ Tate. 
In: 
B.\ Birch, W.\ Kuyk (eds.),
Modular functions of one variable IV,  pp.\ 53--73. 
Springer, Berlin, 1975.
 
\bibitem{Ekedahl 1988}
T.\ Ekedahl:
Canonical models of surfaces of general type in positive characteristic. 
Inst.\ Hautes \'Etudes Sci.\ Publ.\ Math.\ 67 (1988), 97--144.

\bibitem{Fanelli; Schroeer 2020a}
A.\ Fanelli, S.\ Schr\"oer:
Del Pezzo surfaces and Mori fiber spaces in positive characteristic.
Trans.\ Amer.\ Math.\ Soc.\  373 (2020),  1775--1843.

\bibitem{Fanelli; Schroeer 2020b}
A.\ Fanelli, S.\ Schr\"oer:
The maximal unipotent finite quotient, unusual torsion in Fano threefolds, and exceptional Enriques surfaces.
\'Epijournal Geom.\ Alg\'ebrique  4 (2020), Art.\ 11.
 
\bibitem{Grothendieck 1957}
A.\ Grothendieck:
Sur quelques points d'alg\`ebre homologique. 
Tohoku Math.\ J.\ 9 (1957), 119--221. 

\bibitem{EGA IVd}
A.\ Grothendieck:
\'El\'ements de g\'eom\'etrie alg\'ebrique IV: \'Etude locale des
sch\'emas et des morphismes de sch\'emas.
Publ.\ Math., Inst.\ Hautes \'Etud.\ Sci.\   32 (1967).

\bibitem{GB III}
A.\ Grothendieck:
Le groupe de Brauer III.
In: 
J.\ Giraud (ed.) et al.: Dix expos\'es sur la cohomologie des sch\'emas, pp.\ 88--189.
North-Holland, Amsterdam, 1968.

\bibitem{Hartshorne 1994}
R.\ Hartshorne:
Generalised divisors on Gorenstein schemes.
K-Theory 8 (1994), 287--339.

\bibitem{Hogadi 2009}
A.\ Hogadi:
Products of Brauer--Severi surfaces. 
Proc.\ Amer.\ Math.\ Soc.\ 137 (2009),  45--50. 
 
\bibitem{Katsura 1978}
T.\ Katsura:
On Kummer surfaces in characteristic $2$.
In: 
M.\ Nagata (ed.), Proceedings of the international symposium on  algebraic geometry, pp.\ 525--542. 
Kinokuniya Book Store, Tokyo, 1978. 

\bibitem{Kollar 2005}
J.\ Koll\'ar:
Conics in the Grothendieck ring. 
Adv.\ Math.\ 198 (2005),  27--35.
 
\bibitem{Kondo; Schroeer 2019}
S.\ Kondo, S.\ Schr\"oer:
Kummer surfaces associated with group schemes.
Manuscripta Math.\ 166 (2021), 323--342. 

\bibitem{Lang; Tate 1958}
S.\ Lang, J.\ Tate:
Principal homogeneous spaces over abelian varieties. 
Amer.\ J.\ Math.\ 80 (1958), 659--684. 

\bibitem{Laurent; Schroeer 2024}
B.\ Laurent, S.\ Schr\"oer:
Para-abelian varieties and  Albanese maps.
Bull.\ Braz.\ Math.\ Soc.\ 55 (2024), 1--39.
 
\bibitem{Lorenzini; Schroeer 2020}
D.\ Lorenzini, S.\ Schr\"oer:
Moderately ramified actions in positive characteristic.
Math.\ Z.\ 295 (2020), 1095--1142.  

\bibitem{Lorenzini; Schroeer 2023}
D.\ Lorenzini, S.\ Schr\"oer:
Discriminant groups of wild cyclic quotient singularities. 
Algebra  Number Theory 17-5 (2023), 1017--1068.

\bibitem{Matsumura; Oort 1967}
H.\ Matsumura, F.\ Oort:
Representability of group functors, and automorphisms of algebraic schemes.
Invent.\ Math.\ 4 (1967--68), 1--25. 

\bibitem{Mumford 1966}
D.\ Mumford:
Lectures on curves on an algebraic surface.
Princeton University Press, Princeton, 1966.

\bibitem{Mumford 1970}
D.\ Mumford:
Abelian varieties.
Tata Institute of Fundamental Research Studies in Mathematics 5.
Oxford University Press,  London, 1970.

\bibitem{Mumford 1971}
D.\ Mumford:
Theta characteristics of an algebraic curve.
Ann.\ Sci.\ \'Ecole Norm.\ Sup.\   4 (1971), 181--192. 

\bibitem{Raynaud 1970}
M.\ Raynaud:
Sp\'ecialisation du foncteur de Picard.
Publ.\ Math., Inst.\ Hautes \'Etud.\ Sci.\ 38 (1970), 27--76.

\bibitem{Roessler; Schroeer 2022}
D.\ R\"ossler, S.\ Schr\"oer:
Moret-Bailly families and non-liftable schemes.
Algebr.\ Geom.\ 9 (2022),  93--121. 

\bibitem{Saito 2012}
T.\ Saito:
The discriminant and the determinant of a hypersurface of even dimension.  
Math.\ Res.\ Lett.\ 19 (2012),  855--871. 

\bibitem{Schilson 2018}
B.\ Schilson:
Singularit\"aten von Kummer-Variet\"aten in beliebiger Charakteristik.
Dissertation,  D\"usseldorf (2018),
\url{https://nbn-resolving.org/urn/resolver.pl?urn=urn:nbn:de:hbz:061-20181108-114448-1}.

\bibitem{Schilson 2020}
B.\ Schilson:
Wild singularities of Kummer varieties.  
J.\ Singul.\ 20 (2020), 274--288. 

\bibitem{Schroeer 2007}
S. Schr\"oer:
Kummer surfaces for the selfproduct of the cuspidal rational curve.
J.\ Algebraic Geom.\ 16 (2007), 305--346.

\bibitem{Schroeer 2009}
S.\ Schr\"oer:
The Hilbert scheme of points for supersingular abelian surfaces. 
Arkiv Mat.\ 47  (2009), 143--181.

\bibitem{Schroeer 2017}
S.\ Schr\"oer:
Enriques surfaces with normal K3-like coverings.
J.\ Math.\ Soc.\ Japan. 73 (2021), 433--496.
 
\bibitem{Schroeer 2023}
S.\ Schr\"oer:
There is no Enriques surface over the integers.
Ann.\ of Math.\ 197 (2023), 1--63.

\bibitem{Schroeer 2024}
S.\ Schr\"oer:
Albanese maps for open algebraic spaces.
\arXiv{arXiv:2204.02613}, to appear in Int.\ Math.\ Res.\ Not.\ IMRN.

\bibitem{Schroeer; Tziolas 2024}
S.\ Schr\"oer, N.\ Tziolas:
The structure of Frobenius kernels for automorphism group schemes.
Algebra Number Theory 17 (2023),  1637--1680. 

\bibitem{Shioda 1974}
T.\ Shioda:
Kummer surfaces in characteristic $2$. 
Proc.\ Japan Acad.\ 50 (1974), 718--722.

 
\end{thebibliography}
\end{document}